\documentclass[referee,sn-mathphys-num]{sn-jnl}

\usepackage{graphicx}%
\usepackage{multirow}%
\usepackage{amsmath,amssymb,amsfonts}%
\usepackage{amsthm}%
\usepackage{mathrsfs}%
\usepackage[title]{appendix}%
\usepackage{xcolor}%
\usepackage{textcomp}%
\usepackage{manyfoot}%
\usepackage{booktabs}%
\usepackage{algorithm}%
\usepackage{algorithmicx}%
\usepackage{algpseudocode}%
\usepackage{listings}%
\usepackage{longtable}%
\usepackage{setspace}%
\usepackage{float}%
\usepackage{amsthm,mathtools,array}
\usepackage{caption}
\usepackage{subcaption}
\usepackage{tabularx}
\usepackage{multirow}

\theoremstyle{thmstyleone}%
\theoremstyle{thmstyletwo}%

\theoremstyle{thmstylethree}%

\begin{document}

\title[Optimizing Vaccine Site Locations While Considering Travel Inconvenience and Public Health Outcomes]{Optimizing Vaccine Site Locations While Considering Travel Inconvenience and Public Health Outcomes}

\author*[1]{\fnm{Suyanpeng} \sur{Zhang}}\email{suyanpen@usc.edu}

\author[1]{\fnm{Sze-Chuan} \sur{Suen}}\email{ssuen@usc.edu}

\author[1]{\fnm{Han} \sur{Yu}}\email{hyu376@usc.edu}

\author[1]{\fnm{Maged} \sur{Dessouky}}\email{maged@usc.edu}

\author[2]{\fnm{Fernando} \sur{Ordonez}}\email{fordon@uchile.cl}

\affil*[1]{\orgdiv{Daniel J. Epstein Department of Industrial and Systems Engineering, Viterbi School of Engineering}, \orgname{University of Southern California}, \orgaddress{\street{3715 McClintock Ave}, \city{Los Angeles}, \postcode{90089}, \state{CA}, \country{USA}}}

\affil[2]{\orgdiv{Department of Industrial Engineering}, \orgname{Universidad de Chile}, \orgaddress{\street{Beauchef 850}, \city{Santiago}, \postcode{8370451}, \state{Región Metropolitana}, \country{Chile}}}




 \abstract{During the COVID-19 pandemic there were over three million infections in Los Angeles County (LAC). To facilitate distribution when vaccines first became available, LAC set up six \textit{mega-sites} for dispensing a large number of vaccines to the public. To understand if another choice of \textit{mega-site} location would have improved accessibility and health outcomes, and to provide insight into future vaccine allocation problems, we propose a multi-objective mixed integer linear programming model that balances travel convenience, infection reduction, and equitable distribution. We provide a tractable objective formulation that effectively proxies real-world public health goals of reducing infections while considering travel inconvenience and equitable distribution of resources. Compared with the solution empirically used in LAC in 2020, we recommend more dispersed \textit{mega-site} locations that result in a 26\% reduction in travel inconvenience and avert an additional 200 infections.}

\keywords{Vaccine sites, Vaccine allocation, COVID-19, Travel inconvenience}



\maketitle
\section*{Highlights}
\begin{itemize}
    \item We formulate an optimization model to select vaccination sites that consider commuting patterns and public health outcomes.
    \item We devise a tractable expression for incorporating health outcomes into the optimization objective that takes into account vaccination targets, prioritization over areas with a higher risk of becoming infectious, and equitable distribution.
    \item  We implement the optimization model with empirical COVID-19 and traffic data in Los Angeles County. We find that with a proper setup, both travel inconvenience and public health outcomes can be improved from the empirically used solution.

\end{itemize}
\newpage
\section{Introduction}
The COVID-19 pandemic resulted in over three million infections in Los Angeles County (LAC), imposing significant burdens on both public health and the economy for the entire population \cite{richards2022economic}. Months after the COVID-19 outbreak, vaccines were developed to control the spread of COVID-19 as well as prevent severe illness \cite{alagoz2021impact,huang2022vaccination}.

When COVID-19 vaccines first became available and distribution to local pharmacies had not yet been established, the LAC Department of Public Health needed to quickly and efficiently distribute a large volume of vaccines to a population of 10 million LAC residents. LAC therefore initiated a vaccination program which used large venues such as parking lots, universities, and theme parks, referred to as `mega-sites,' to facilitate mass vaccination efforts \cite{megasiteslac}. LAC used a total of six \textit{mega-sites} between January to April, 2021. LAC had a sign-up system where residents provided demographic information and signed up for available \textit{mega-site} appointments \cite{myturn}; availability at different sites and times could be restricted depending on demographic characteristics such as age and work industry. This system essentially assigned residents to mega-sites according to rules which changed over time as vaccination restrictions were lifted to prioritize vulnerable residents \cite{vaccinepriority}. In this period, each \textit{mega-site} was able to administer between 200,000 to 800,000 vaccinations per month. 

The locations of \textit{mega-sites} should be strategically chosen to encourage vaccination for populations in high-risk infection areas, thereby reducing the disease's spread. Additionally, these sites should be placed to reduce the overall travel inconvenience for individuals, since these additional time costs could present an additional barrier to vaccine access. 

Prior work has shown that shorter trips that require less travel inconvenience and less time off work may lead to greater vaccine access \cite{papastergiou2014community}. However, the distribution of \textit{mega-sites} in LAC was generally concentrated near downtown, resulting in long travel times for suburban residents. For example, residents in the northeast corner of LAC would need to travel as far as 60 miles to reach the nearest \textit{mega-sites}. Reducing travel inconvenience to vaccination centers may expedite the uptake of vaccines, which is critical for controlling outbreaks \cite{papastergiou2014community,lacombe2018human,schmitzberger2022identifying}. Strategically placed \textit{mega-sites} may therefore improve vaccination coverage. 

However, prior work on healthcare facility location and health transportation typically focus on transportation costs, where transportation costs are typically measured as the distance from the home to the medical service center. This generally does not consider patients' commuting patterns, which may mean this travel inconvenience in traditional literature may not reflect actual inconvenience to vaccine sites. 

Exploiting commuting patterns to place health facilities may improve both problem realism and allow for better transportation solutions. Identifying realistic optimal locations for vaccination requires considering how vaccination would be integrated into a patient's daily life. In particular, travel inconveniences can occur in various situations for commuters, such as when commuting to a vaccination site while on the way to or from work, or when traveling from work or home to a vaccine site and then returning to work or home afterward. LAC residents are notorious for regularly spending long periods of time traveling for work and entertainment; in fact, empirical data shows that workers in LAC travel an average of 54 minutes daily from their residence to their workplace \cite{surveylacommute}. Commuting patterns must therefore be integrated into decisions for optimal placement of \textit{mega-site} locations, as individuals may travel to vaccine sites before or after work, and vaccine sites can be conveniently placed next to workplaces or close to home. 

In this work, we therefore consider travel inconvenience, which considers commuting patterns, instead of total transportation costs from the home only as in traditional facility location problems. Considering commuting patterns in such an optimization problem is particularly timely, as large datasets on traffic volume, time, and speed have already been collected using road sensor data since 2013 in LAC \cite{pems}. While this data has primarily been used for traffic congestion and transportation research, it can be leveraged for critical insights into health resource distribution applications. 

It is also important to consider health outcomes in addition to travel inconvenience when choosing vaccine \textit{mega-site} locations. Residents from a high infection risk area should be able to be assigned to a convenient vaccination \textit{mega-site} with low travel inconvenience to improve health outcomes such as the number of infections averted \cite{haas2022infections}. In this paper, we study how to optimize the locations of \textit{mega-sites} and vaccine assignment for minimizing travel inconvenience and achieving desirable public health outcomes.

To do this, it is critical to understand which group should be prioritized for vaccination in order to mitigate infections. Prior work has demonstrated that geographically targeted prioritization vaccination strategies for vaccine allocation to groups at high risk of becoming infected can be effective at reducing disease burden  \cite{araz2012geographic,wrigley2021geographically}. Measuring social interaction levels can be one measure to identify groups at a higher risk of infection \cite{karmakar2021association}. Similar to \cite{hu2021human}, we use regional mobility and traffic flow patterns as a proxy measure for social interaction level to identifying regional subpopulations at higher risk of infection. 

It is therefore important to consider commuting patterns and regional vaccine allocation prioritization when choosing vaccination \textit{mega-sites} to promote vaccine accessibility and reduce disease cases. While many works have studied vaccine allocation and facility location problems for infectious disease, we here take a novel approach by integrating travel inconvenience, public health outcomes, and preferences for equitable accessibility into an optimization model to assist public health decisions around determining vaccine supply locations. 

\subsection{Literature Review}
Prior works have studied vaccination allocation problems as well as locations of vaccination sites. However, to the best of our knowledge, they generally do not integrate both travel inconvenience and public health outcomes into the optimization objective. We review several relevant works here and highlight how our approach differs from these existing works.  

There is a rich literature on facility location problems in healthcare. Traditional research on the location of healthcare facilities has primarily focused on accessibility and coverage \cite{daskin2004location,ahmadi2017survey, marianov2001optimal,gunecs2019location}. Additionally, certain studies on coverage exclude travel inconvenience, focusing instead on probabilistic demand distributions \cite{jia2007modeling,silva2008locating,chan2016optimizing}. Besides coverage problems, \cite{mete2010stochastic,zhang2015novel} have developed models for locating medical service stations for emergency situations which focus on estimating demand and routing medical service units while considering inventory control. This contrasts with our approach where individuals travel to a vaccination center with a constant set capacity over time. \cite{beheshtifar2015multiobjective} provided a multiobjective optimization formulation for solving the location and resource allocation problem with a focus on both transportation and equitable access. In this work, we extend this idea to include disease dynamics, expanding the traditional location problems to include multiple objectives on disease outcomes, travel inconvenience, and equitable access. For a more detailed review of works on location problems in healthcare applications, we refer the reader to \cite{ahmadi2017survey}. In general, prior studies did not integrate commuting patterns in optimizing the locations of vaccine sites \cite{rastegar2021inventory,chen2022location,valizadeh2023designing}, which we do here. There is also a wide range of general facility location problems that we are unable to cover comprehensively in this review. These models often address various factors beyond travel inconvenience but are not directly applicable to healthcare facility location problems. For instance, the placement of electric vehicle charging stations considers both travel time and congestion costs \cite{ahmad2022optimal,huang2020electric,kchaou2021charging}. Similarly, retail planning focuses more on meeting customer demand rather than healthcare-specific considerations \cite{drezner2002retail,muller2014customer,hirpara2021retail}. For a more comprehensive review of facility location problems, we refer readers to \cite{celik2020comparative}.

A separate literature has considered health outcomes in vaccination distribution problems, but generally focuses on resource allocation as opposed to locating distribution centers \cite{buckner2020optimal,choi2020optimal,deng2021joint,mukandavire2020quantifying,foy2021comparing,shim2021optimal,becker1997optimal,aaby2006montgomery,bertsimas2020optimizing,buhat2021using}. For example, \cite{bertsimas2020optimizing} optimizes the allocation of COVID-19 vaccines by state while incorporating prevalence predictions on the future epidemic. \cite{buhat2021using} developed a linear programming formulation for optimizing vaccine allocation in the Philippines. A prioritization score is assigned to each region and is computed based on the local basic reproduction number. These studies concentrate on the allocation of resources and their impact on public health outcomes. However, they do not take into account travel inconvenience, as the placement of vaccination sites is either predetermined or not a focus of their analyses. In this paper, we not only incorporate commuting patterns into the calculation of travel inconvenience for determining vaccination site locations, but we also derive a straightforward formula that effectively captures the dynamics of infectious diseases resulting in a more efficient solution for addressing the vaccine distribution challenge while taking into account the complexities of disease spread.

We consider regional vaccine allocation prioritization in our formulation. Geographical prioritization strategies for vaccination allocation have demonstrated their efficiency in safeguarding groups at high risk of becoming infected \cite{araz2012geographic,wrigley2021geographically,karmakar2021association,wilder2020allocating,bubar2021model,choi2021vaccination,tatapudi2021impact,lee2015vaccine,lee2021strategies}. Numerous techniques have been proposed to geographically prioritize groups at high risk of infection \cite{wilder2020allocating,chang2021mobility,zachreson2021risk}. In this paper, we use traffic flow to proxy for infection risk from increased social contacts as the mobility patterns are strongly correlated with decreased COVID-19 case growth rates \cite{badr2020association}, and we prioritize vaccination accessibility for regions that have greater transportation activity. 

\subsection{Contributions}
We make several contributions in this study. We formulate an optimization model to select vaccination sites that consider commuting patterns and public health outcomes. While substantial prior literature studies similar facility location problems, we are not aware of any prior work that has considered commuting patterns in such an optimization model. Also, we devise a tractable expression for incorporating health outcomes into the optimization objective that takes into account vaccination targets, prioritization over areas with a higher risk of becoming infectious, and equitable distribution. Specifically, we propose a way to compute vaccination uptake targets (to reach herd immunity) for different geographical areas. Then we evaluate several different vaccine prioritization strategies and identify the one that is most able to reduce infections. We additionally incorporate preferences around the equitable distribution of vaccinations into the objective function to ensure smoothness in vaccination uptake across different regions. Moreover, we implement the optimization model with empirical COVID-19 and traffic data in LAC, the most populous county in the U.S. and an epicenter in the COVID-19 pandemic, and we compare outcomes using the optimal \textit{mega-site} locations with those used in LAC 2020. We demonstrate that with a proper setup, both travel inconvenience and public health outcomes can be improved from the empirically used solution. This work provides insights to policy-makers when setting up large vaccine sites for controlling transmissible diseases at the initial phase of vaccination in the future.

The remainder of this paper is organized as follows: we present the problem setup and the optimization formulation in Section~\ref{formulation}. The numerical example is shown in Section~\ref{numeric}. In Section~\ref{conclusion}, we conclude.

\section{Model Formulation}\label{formulation}
In this section, we present an optimization formulation for optimizing the location of \textit{mega-sites} and allocating vaccines to the public. We assume a disease outbreak setting and vaccine is newly introduced with limited supply. In the optimization problem, we focus on two questions: (1) where should the public health department set up \textit{mega-sites} for large vaccine dispensing; and (2) how to assign individuals for vaccination. We assume that the \textit{mega-site} locations will not change during the analysis period. We assume each \textit{mega-site} has a maximum number of vaccines to allocate in each time epoch. We assume each region can have at most one \textit{mega-site}. Regions can be defined in a variety of ways, and one convenient way may be using pre-existing definitions by health decision-makers; for instance, in LAC, `Health Districts' are used to delineate health resources by geography \cite{guevara2015changing}. We will therefore use the term Health District (HD) to represent regions within an area. 

\subsection{Minimization of Travel Inconvenience Considering Commuting Patterns}\label{transport_model} A model that does not account for commuting patterns may be imprecise in estimating the travel inconvenience. Therefore, in this section, we consider two categories of individuals in our population for vaccine assignments: commuters and non-commuters. This allows a more precise estimation of cost quantification. In this study, we use HD to represent locations and let $\mathcal{H}$ denote the set of all HDs. We use $e_{uv}$ to denote the number of commuters traveling from HD $u \in \mathcal{H}$ to HD $v \in \mathcal{H}$. The non-commuter population in each HD $u \in \mathcal{H}$ is defined as the total population, denoted by $p_u$, minus the commuters from HD $u$. We compute travel inconvenience for commuters based on distances from vaccine sites to both their residence and their workplaces.  For instance, if an individual resides in HD $u$ and commutes to their workplace in HD $v$, a vaccination site situated near either HD $u$ or HD $v$ should be regarded as convenient for that person. 

 The additional travel inconvenience for non-commuters is $c_{uw}+c_{wu}$ for people who live in HD $u\in\mathcal{H}$ and get vaccinated at HD $w\in \mathcal{H}$. For commuters, four different possible cases can be considered. An individual who lives in HD $u\in \mathcal{H}$, goes to work in HD $v\in \mathcal{H}$, and gets vaccinated at HD $w\in \mathcal{H}$ can choose to get vaccinated either before or after work, during work hours, en route to work, or while returning from work. We assume they will always choose the option that minimizes their total travel inconvenience. We use $d_{uvw}$ to denote the additional travel inconvenience (minutes) for people who live in HD $u\in \mathcal{H}$, go to work in HD $v\in \mathcal{H}$, and get vaccinated at HD $w\in \mathcal{H}$. Therefore, the additional travel inconvenience for this individual is defined as $d_{uvw}=\min\{c_{uw}+c_{wu},c_{vw}+c_{wv}, c_{vw}+c_{wu}-c_{vu},c_{uw}+c_{wv}-c_{uv}\}$, which is known for each commuter. 

 We consider three decision variables for each time period within the set of all time periods $\mathcal{T}$. Let $x_u$ represent the decision on whether to open a \textit{mega-site} in HD $u \in \mathcal{H}$. The variable $y_{uw}^t$ denotes the number of non-commuters residing in HD $u \in \mathcal{H}$ that are assigned to get vaccinated at HD $w \in \mathcal{H}$ during decision period $t \in \mathcal{T}$. Similarly, $z_{uvw}^t$ represents the number of commuters who live in HD $u \in \mathcal{H}$, work in HD $v \in \mathcal{H}$, and are assigned to get vaccinated at HD $w \in \mathcal{H}$ during decision period $t \in \mathcal{T}$.

Furthermore, the total number of \textit{mega-sites} is restricted to not exceed $K$. Finally, let $U_{tw}$ denote the number of vaccines available at \textit{mega-site} $w \in \mathcal{H}$ in decision period $t$. The detailed description of the notation is described in Table~\ref{tab:table_of_notations}.

\begin{table}[h]
\setlength\extrarowheight{2pt}
\begin{tabularx}{\textwidth}{ |m{2cm}|c|X| } 
\hline
 & Notation & Definition \\
\hline
\multirow{2}{4em}{Set} &$\mathcal{H}$& a set of HDs. \\ 
& $\mathcal{T}$& a set of decision periods. \\ \hline
\multirow{10}{4em}{Parameters} & $K$& total number of \textit{mega-sites} that can be opened. \\ 
&$c_{uv}$& the travel inconvenience (minutes) from HD $u\in \mathcal{H}$ to HD $v\in \mathcal{H}$.\\
&  $d_{uvw}$& the additional travel inconvenience (minutes) for people who live in HD $u\in \mathcal{H}$, go to work in HD $v\in \mathcal{H}$, and get vaccinated at HD $w\in \mathcal{H}$. \\ 
& $p_u$& total population at HD $u\in \mathcal{H}$. \\ 
& $e_{uv}$& number of people who live in HD $u\in \mathcal{H}$ and go to work in at HD $v\in \mathcal{H}$.   \\ 
& $U_{tw}$& number of vaccines available at the \textit{mega-site} $w\in \mathcal{H}$ for decision period $t \in \mathcal{T}$. \\ \hline
\multirow{6}{4em}{Variables} & $x_u$& open a \textit{mega-site} in HD $u\in\mathcal{H}$. \\ 
& $y_{uw}^t$& number of non-commuters who live in HD $u\in \mathcal{H}$ that are assigned to get vaccinated at HD $w \in \mathcal{H}$ during the decision period $t\in \mathcal{T}$.\\ 
& $z_{uvw}^t$ &  number of commuters who live in HD $u\in \mathcal{H}$, go to work in HD $v\in \mathcal{H}$ that are assigned to get vaccinated at HD $w \in \mathcal{H}$ during the decision period $t\in \mathcal{T}$.\\
\hline
\end{tabularx}
\caption{Table of notation.}
\label{tab:table_of_notations}
\end{table}

We consider a traditional mixed integer linear programming (MILP) problem that minimizes additional travel inconvenience:
\begin{subequations}
    \label{P1_formulation} 
    \begin{equation*}
  \text{(P1)}\quad  \min_{x,y,z} \quad \sum_{t\in \mathcal{T}}\sum_{u\in \mathcal{H}}\sum_{w\in \mathcal{H}}y_{uw}^t(c_{uw}+c_{wu})+ \sum_{t\in \mathcal{T}}\sum_{u\in \mathcal{H}}\sum_{v\in \mathcal{H}}\sum_{w\in \mathcal{H}}z_{uvw}^td_{uvw} 
  \end{equation*}
  \begin{align}
    s.t. \quad& \sum_{u\in \mathcal{H}}x_u \leq K\\
    &\sum_{t\in \mathcal{T}}\sum_{w\in \mathcal{H}}y_{uw}^t = p_u-\sum_{v\in \mathcal{H}}e_{uv},\forall u \in \mathcal{H}\\
    & \sum_{t\in \mathcal{T}}\sum_{w\in \mathcal{H}}z_{uvw}^t=e_{uv},\forall u\in \mathcal{H},v\in \mathcal{H}\\
    & y_{uw}^t\leq U_{tw}x_{w},\forall u\in \mathcal{H},w\in \mathcal{H},t \in \mathcal{T}\\
    & z_{uvw}^t\leq U_{tw}x_{w},\forall u\in \mathcal{H},v\in\mathcal{H},w\in \mathcal{H},t \in \mathcal{T}\\
    & \sum_{u\in \mathcal{H}}y_{uw}^t+\sum_{u\in \mathcal{H}}\sum_{v\in \mathcal{H}}z_{uvw}^t\leq U_{tw},\forall w\in \mathcal{H},t\in \mathcal{T}\\
    & x_u \in \mathbb{B},\forall u\in \mathcal{H}\\
    & y_{uw}^t\in \mathbb{Z}^+,\forall u\in \mathcal{H},w\in \mathcal{H},t \in \mathcal{T}\\
    & z_{uvw}^t\in \mathbb{Z}^+,\forall u\in \mathcal{H},v\in \mathcal{H},w\in\mathcal{H},t \in \mathcal{T}
    \end{align}
\end{subequations}

The objective in P\ref{P1_formulation} is minimizing the total travel inconvenience from home/workplaces to vaccine sites. Constraint 1a limits the number of \textit{mega-sites} that can be built. Constraints 1b and 1c ensure that all individuals are vaccinated eventually. In this paper, we assume that all people eventually get vaccinated in one of the \textit{mega-sites} during one of the decision periods $t \in \mathcal{T}$. Constraints 1d and 1e ensure that individuals can only get vaccinated in an HD where a \textit{mega-site} exists. Constraint 1f restricts the number of vaccines that can be administrated by the \textit{mega-site} at the HD $w\in\mathcal{H}$ during the decision period $t\in \mathcal{T}$. Constraints 1d and 1e are formulated using big-M notation. We implement this model in Gurobi 9.5.2 using Python 3.9 and describe the outcomes in Section~\ref{numeric}.

We consider the traditional transportation model that does not include commuting patterns (treat all individuals as non-commuters, i.e., $e_{uv}=0,\forall u,v\in\mathcal{H}$) as model P0. We compare P\ref{P1_formulation} with P0 in the numerical section.

\subsection{Modeling Disease Dynamics}
The MILP formulation in P\ref{P1_formulation} only minimizes the total travel inconvenience. However, we are also interested in reducing poor health outcomes due to infectious disease in this population. In the next section, we therefore integrate the public health objective into the optimization model. 

\subsubsection{Compartmental Model}\label{compart_intro}

Vaccines can prevent downstream transmission of disease, so it is important to track disease dynamics over time in order to understand the impact of \textit{mega-site} location choices and allocation strategies. Compartmental models are commonly used for modeling complex disease dynamics of infectious disease on a population level \cite{buckner2020optimal,choi2020optimal,deng2021joint,mukandavire2020quantifying,foy2021comparing,shim2021optimal}. In this paper, we formulate a compartmental model that tracks the disease dynamics in different HDs that helps us integrate public health objectives into the optimization model. Due to the complexity of disease dynamics, compartmental models are hard to formulate and often intractable except through simulation. We therefore use a simple susceptible-infectious-recovered (SIR) model with vaccinations which simplifies the disease dynamics to create a tractable objective term, we then use a more detailed compartmental model for evaluating the performance of our model. SIR models have been used extensively in many disease contexts and track health states, transmission, and recovery over time \cite{weiss2013sir,cooper2020sir,tang2020review}. Our model is additionally stratified by HD to track the number of individuals within each HD in each health state.

In the model, we use $s^v,i^v,r^v$ to represent the proportion of individuals in each health state (susceptible ($s$), infectious ($i$), and recovered ($r$)) in HD $v\in\mathcal{H}$. We capture the transmission rate between HDs using $\beta_{uv},\forall u\in\mathcal{H},\forall v \in \mathcal{H}$, which denotes the transmission rate from HD $u$ to HD $v$. We assume $\gamma$ is the geographically invariant clearance rate. Our compartmental model is then:
\begin{align*}
    &\frac{ds^{v}}{dt}=- \sum_u {\beta}_{uv} s^v (t) {{i}^u(t)}\\
    &\frac{di^v}{dt} =  \sum_u {\beta}_{uv} s^v (t){{i}^u(t)}-{\gamma i^v(t)}\\
    &\frac{dr^v}{dt} = \gamma i^v(t)\\
\end{align*}

\subsubsection{Herd Immunity Threshold Analysis}\label{method:herd_immunity}

Much of the benefits of a population vaccination campaign can be garnered before all individuals have been vaccinated. In epidemiology, the \textit{herd immunity threshold} is the idea that less than 100\% of the population need be immune to the disease to prevent an epidemic \cite{chowdhury2021universality}. We will use this concept to formulate our vaccine threshold targets. 

Many works have studied computation methods for identifying herd immunity thresholds \cite{fine2011herd,aguas2020herd,chowdhury2021universality} within a population. To compute the herd immunity threshold in this paper, we use a similar approach as \cite{chowdhury2021universality}. We find the herd immunity threshold such that the spread of disease will start to decrease ($\frac{di^v}{dt}<0, \forall v\in \mathcal{H}$). 

Then we can write the condition for decreasing prevalence as: 
\begin{equation*}
    \frac{di^v}{dt} = \sum_{u}\beta_{uv} s^v(t)i^u(t) -\gamma i^v(t)<0,\forall v \in \mathcal{H},
\end{equation*}
where $\beta_{uv}$ is the transmission rate and $\gamma$ is the recovery rate.

This is achieved when  
\begin{equation*}
s^v(t)<\frac{\gamma i^v(t)}{\sum_{u}\beta_{uv} i^u(t)},\forall v \in \mathcal{H}.
\end{equation*}

The herd immunity threshold is the lowest herd immunity level required for the epidemic to decrease. Although herd immunity can be achieved through a rise in the number of recovered or via vaccinations \cite{chowdhury2021universality}, in our paper, we presume that herd immunity is primarily attained through vaccination, given that the number of infected and recovered individuals are relatively small compared to the entire population, particularly as we are interested in the initial stages of a pandemic. Let us define $L^v$ to be the herd immunity level that leads the epidemic to decrease. Then $L^v$ can be approximated by $1-s^v(t)$ at time $t$ for region $v$. 
\begin{equation*}
    L^v=1-s^v(t)>1-\frac{\gamma i^v(t)}{\sum_{u}\beta_{uv} i^u(t)},\forall v \in \mathcal{H}.
\end{equation*}

 then simplify by assuming that the ratio of prevalences (denoted as $\kappa$) between HDs is a time-invariant constant and this has been validated in the numeric simulation. We use the average ratio over a specific period $[t^-, t^+]$ to represent this constant. Specifically, we define:

$$\kappa_{uv} = \frac{1}{\text{\# of time points in }[t^- ,t^+]} \sum_{t=t^-}^{t^+} \frac{i^u(t)}{i^v(t)}, \quad \forall u \in \mathcal{H}, \forall v \in \mathcal{H},$$
where time points within the interval $[t^-, t^+]$ are uniformly discretized.

Then the previous equation can be rewritten as 
\begin{equation*}
     L^v>\frac{\sum_u \beta_{uv}\kappa_{uv}-\gamma}{\sum_u \beta_{uv}\kappa_{uv}},\forall v\in \mathcal{H}.
\end{equation*}

We therefore set the vaccination target for each HD $u$ as the herd immunity threshold $L^{u*}=\frac{\sum_{v\in\mathcal{H}} \beta_{vu}\kappa_{vu}-\gamma}{\sum_{v\in\mathcal{H}} \beta_{vu}\kappa_{vu}}$. We incorporate these vaccination targets $L^{u*}$ into our optimization objective. We wish to minimize the difference between vaccination levels and $L^{u*}$ targets across all time: $p_u L^{u*}-\sum_{t'=1}^{t}(\sum_w y_{uw}^{t'}+\sum_v\sum_w z_{uvw}^{t'}), u\in\mathcal{H},t\in \mathcal{T}$. We add this term to the optimization problem to make it a multi-objective formulation. 

This objective considers equal weights for each HD in the allocation process. However, we may want to prioritize HDs with residents at higher risk of becoming infected; in this work, we will consider transmission risk proxied by greater transportation movement around the city. In general, weighting each HD by a risk score could capture a wide variety of definitions of risk (due to demographic characteristics, e.g., age, or differences in access to treatment, etc.). We term these weights `prioritization scores' as they allow the optimization model to prioritize specific regions and allocate vaccinations to these areas preferentially. Such prioritization scores for areas with different risks have been used in prior literature to help capture heterogeneity in disease control \cite{wilder2020allocating,chen2021prioritizing}. 

We use $\delta_u^t$ to define the prioritization score for HD $u\in\mathcal{H}$ at time period $t\in\mathcal{T}$. We let our prioritization score $\delta_u^t$ be the flow-out in HD $u$ divided by the total flow across all HDs. In Section~\ref{numeric}, we examine the performance of our prioritization score and compare it against other common scores in the literature.

\subsubsection{Equitable Distribution of Vaccines}
We also need to consider vaccination distribution equity. Access to healthcare resources without distinction between individuals is a fundamental tenet of promoting human rights \cite{aborode2021equal,bollyky2020equitable}. Thus, in our optimization model, we aim to reduce the disparity in the number of individuals vaccinated across various HDs during the initial vaccination phase ($\mathcal{T}_E\subset \mathcal{T}$) when supply is mostly available through \textit{mega-sites} and not available in local clinics. 
As the proportion of vaccinated individuals grows and the urgency of vaccination wanes due to protective effects from herd immunity, equity may be a less pressing concern. Of course, $\mathcal{T}_E$ can be set equal to $\mathcal{T}$ if equity objectives should be considered throughout the time horizon. 


\subsection{An Multi-Objective Model} \label{multi-obj}
Realizing the need to prioritize high risk groups of becoming ill, promote equitable distribution of vaccines, as well as improve accessibility to vaccination, we formulate another model that incorporates these objectives. Let $\zeta^t_u$ and $\tau_t$ present the objective value on public health outcome at time $t$ for HD $u$ and equitable distribution at time $t$ respectively, we have a MILP formulation as follows:
\begin{subequations}
\label{P2_formulation}
\begin{equation*}
   \text{(P2)} \quad \min_{x,y,z} \quad \sum_{t\in \mathcal{T}}\sum_{u\in \mathcal{H}}\sum_{w\in \mathcal{H}}y_{uw}^t(c_{uw}+c_{wu})+ \sum_{t\in \mathcal{T}}\sum_{u\in \mathcal{H}}\sum_{v\in \mathcal{H}}\sum_{w\in \mathcal{H}}z_{uvw}^td_{uvw}+\lambda \sum_{t\in\mathcal{T}} \sum_{u\in\mathcal{H}} \delta^t_u \zeta_u^t+\lambda' \sum_{t\in \mathcal{T}_E}\tau_t 
   \end{equation*}
   \begin{align}
    s.t. \quad& \sum_{u\in \mathcal{H}}x_u \leq K\\
    &\sum_{t\in \mathcal{T}}\sum_{w\in \mathcal{H}}y_{uw}^t = p_u-\sum_{v\in \mathcal{H}}e_{uv},\forall u \in \mathcal{H}\\
    & \sum_{t\in \mathcal{T}}\sum_{w\in \mathcal{H}}z_{uvw}^t=e_{uv},\forall u\in \mathcal{H},v\in \mathcal{H}\\
    & y_{uw}^t\leq U_{tw}x_{w},\forall u\in \mathcal{H},w\in \mathcal{H},t \in \mathcal{T}\\
    & z_{uvw}^t\leq U_{tw}x_{w},\forall u\in \mathcal{H},v\in \mathcal{H},w\in \mathcal{H},t\in \mathcal{T}\\
    & \sum_{u\in \mathcal{H}}y_{uw}^t+\sum_{u\in \mathcal{H}}\sum_{v\in \mathcal{H}} z_{uvw}^t\leq U_{tw},\forall w\in \mathcal{H},t\in \mathcal{T}\\
    & \zeta_u^t\geq p_u L^{u*}-\sum_{t'=1}^t(\sum_{w\in\mathcal{H}} y_{uw}^{t'}+\sum_{v\in\mathcal{H}}\sum_{w\in\mathcal{H}} z_{uvw}^{t'}),\forall u\in\mathcal{H},t\in\mathcal{T}\\
	&\tau_t \geq \frac{(\sum_{v\in\mathcal{H}}\sum_{w\in\mathcal{H}} z_{uvw}^t+\sum_{v\in\mathcal{H}} y_{uv}^t)}{p_u}-\frac{(\sum_{v\in\mathcal{H}}\sum_{w\in\mathcal{H}} z_{u'vw}^t+\sum_{v\in\mathcal{H}} y_{u'v}^t)}{p_{u'}},  \forall u\in\mathcal{H},u'\in\mathcal{H},t\in \mathcal{T}_E\\
    & x_u \in \mathbb{B},\forall u\in \mathcal{H}\\
    & y_{uw}^t\in \mathbb{Z}^+,\forall u\in \mathcal{H},w\in \mathcal{H},t \in \mathcal{T}\\
    &z_{uvw}^t\in \mathbb{Z}^+,\forall u\in \mathcal{H},v\in \mathcal{H},w\in \mathcal{H}, t\in \mathcal{T}\\
    & \zeta_u^t\geq 0,\forall u \in \mathcal{H},t\in\mathcal{T}
    \end{align}
\end{subequations}

Unlike P\ref{P1_formulation} in Section~\ref{transport_model}, this new MILP model (P\ref{P2_formulation}) includes a public health objective as well as an equitable distribution objective. In the objective, $\zeta_u^t$ is the remaining unvaccinated population at HD $u\in\mathcal{H}$ at time $t\in\mathcal{T}$ relative to the vaccination threshold target $L^{u*}$. $\delta_u^t$ denotes the prioritization score for HD $u\in\mathcal{H}$ at time $t\in\mathcal{T}$. $\tau_t$ is the maximal difference in the number of vaccines distributed to two HDs at period $t$ in the early stages of vaccination $\mathcal{T}_E\subset\mathcal{T}$. $\lambda'$ balances the importance of equitable distribution and prioritization over groups with higher risks of becoming infectious. $\lambda$ balances the importance of travel inconvenience and public health outcomes. Constraints 2g and 2h compute the value of $\zeta_u^t,\tau_t$ from the vaccine allocation strategy. 

In Section~\ref{numeric}, we parameterize this optimization problem using data from LAC and solve it using Gurobi 9.5.2 on Python 3.9. We compare it with the solutions from P0, P\ref{P1_formulation} in Section~\ref{transport_model}, and modeled outcomes using the empirically chosen locations in LAC.

\section{Numerical Example}\label{numeric}

Using the disease modeling framework and the optimization model presented in Section~\ref{formulation}, we outline the model's parameterization, evaluate the model's performance by the number of infections averted, travel inconvenience, and equitable distribution of vaccines among the 26 HDs within LAC, and compare results to those using the empirically chosen locations in 2021.

\subsection{Model Inputs}\label{inputs}
We use traffic data from Performance Measurement System (PeMS) Data Source \cite{pems}, Open Route API \cite{openroute}, and population data \cite{population,commuter} to parameterize the optimization model. We assume six \textit{mega-sites} can be built at maximum. We assume each can provide a capacity of 400,000 monthly vaccinations. The optimization model solves monthly decisions across a six months time horizon. In the rest of Section~\ref{inputs}, we provide a comprehensive explanation of model inputs.

\subsubsection{Commuters and Non-commuters}

\begin{figure}[h]%
\centering
\includegraphics[width=1\textwidth]{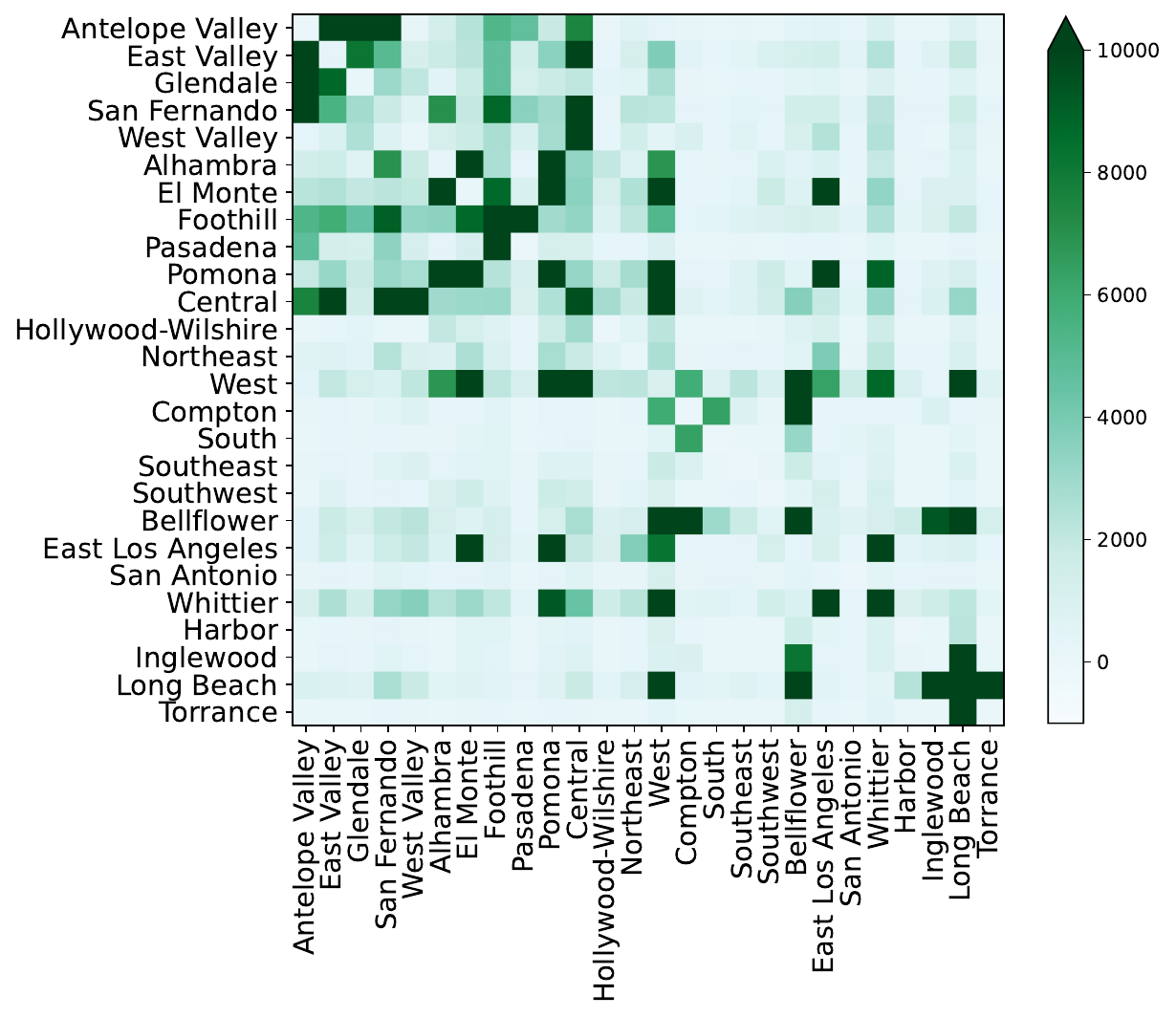}
\caption{Estimation of daily traffic flow using PeMS data. Each entry represents the number of vehicles going from an origin HD (row) to a destination HD (column). A darker color indicates a higher volume of traffic. We assume individuals live in the origin HD and work in the destination HD.}\label{OD_heatmap}
\end{figure}

We use traffic volume data from PeMS database \cite{pems} to compute the number of daily commuters from one HD to another. PeMS contains the number of vehicles passing over a sensor in 5 minutes intervals and the average speed of these vehicles. We approximate the number of vehicles traveling from one region to another region by extending a previously published dynamic origin-destination estimation (DODE) model \cite{ma2018estimating,yu2024extending}. We depict a heatmap of the estimated traffic flow volume by number of vehicles on a daily basis in Figure~\ref{OD_heatmap}. Figure~\ref{OD_heatmap} represents a total of 2,876,523 vehicles. In 2020, there were an estimated 4,396,232 commuters across LAC \cite{commuter}. Given evidence suggesting a dramatic decrease in public transportation usage during that year \cite{liu2020impacts}, we assume that people used cars to commute rather than public transportation. We also assume that people did not walk to mega-sites for vaccination; indeed, LAC offered drive-through-only vaccinations at these empirical mega-sites. Therefore, to obtain a matrix representing the number of commuters in LAC, we scale up the vehicle matrix by multiplying all entries by 1.53—the average number of people per vehicle (calculated by dividing the total number of commuters by the number of vehicles). The number of commuters in each HD can then be computed by summing the corresponding row of the commuter matrix. Then the number of non-commuters can be obtained by taking the difference from the total population in each HD \cite{population}. 

\subsubsection{Parameterization of Travel Cost Matrix}
For an individual who lives in or works in a certain HD, we use the geographic centroid of that HD to represent the home or work address. Similarly, we use the geographic centroid of each HD to approximate the potential location of the \textit{mega-site} in an HD. To parameterize the cost of travel from one HD $u\in\mathcal{H}$ to another HD $v\in\mathcal{H}$, we use Open Route API \cite{openroute} to compute the time (in minutes) needed to traveling from $u$ to $v$ ($c_{uv}$) under traffic-free conditions. We use these values to parameterize the travel cost matrices ($c, d$)  for commuters and non-commuters as described in Section~\ref{transport_model}.

\subsubsection{Herd Immunity Threshold}

We calibrated a compartmental model to estimate the transmission and recovery rates. Using these calibrated rates and the number of infectious individuals from the model, we then calculated the herd immunity threshold. We find the vaccination target for each HD to ensure herd immunity using $L^{u*},\forall u \in \mathcal{H}$ in Section~\ref{method:herd_immunity}. 
High vaccination targets are observed for all HDs due to the relatively high transmission rate compared to the recovery rate. Among all HDs, West Valley has the lowest vaccination target (95.13\%), and Pasadena has the highest vaccination target (99.15\%). On average, 97.36\% of the total population needs to be vaccinated to achieve herd immunity. Detailed herd immunity thresholds across HDs are shown in Appendix~\ref{herd_results}. However, the herd immunity threshold does not dictate prioritization within the population when no health districts (HDs) have reached the threshold or when all HDs have met the target. We therefore introduce the prioritization score in the next section.

\subsubsection{Prioritization Score}\label{score_choice}

In formulation P\ref{P2_formulation}, we weight the herd immunity level needed in each HD by a prioritization score ($\delta_u^t$) to preferentially target more vulnerable HDs first. We prioritize vaccination accessibility for regions that have greater traffic flow. Specifically, we let the score be the number of daily commuters in each HD divided by the total number of commuters in LAC,
\begin{align*}
    \delta_u^t = c(t) \frac{\sum_{v\in \mathcal{H}}e_{uv}}{\sum_{u,v\in \mathcal{H}}e_{uv}} ,
\end{align*}
where $c(t)$ is the time-dependent discount factor. We apply a discount factor to the prioritization score to reduce the value of future vaccination rewards compared to immediate vaccination. In other words, this encourages individuals to get vaccinated as soon as possible. 

In the Appendix, we compare our prioritization score with others from the literature. We found that the score based on commuting patterns resulted in significantly more infections averted than the others, while its performance on the two other metrics was similar. Therefore, we use traffic flow to define the prioritization score.

To sum up, the determination of \textit{mega-sites} requires knowing traffic flow, commuter and non-commuter population sizes, and herd immunity thresholds. Traffic flow is commonly estimated using transportation databases, population statistics are usually obtained from census databases, and herd immunity thresholds are calculated based on data from disease surveillance systems.

\subsection{Evaluating Public Health Outcomes} 

\begin{figure}[h]%
\centering
\includegraphics[width=0.5\textwidth]{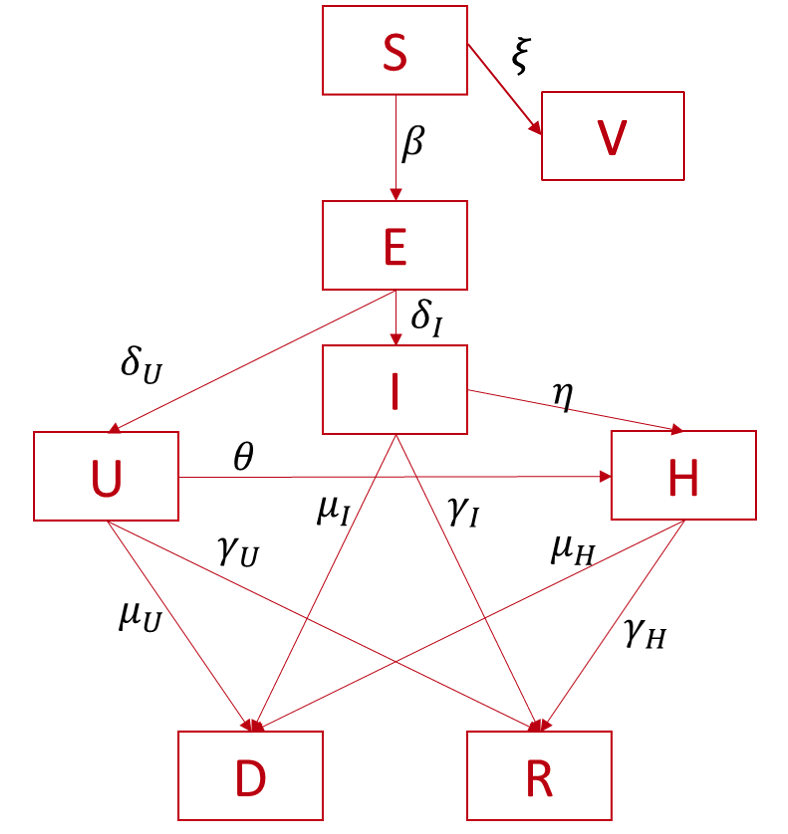}
\caption{The compartmental model for modeling COVID-19 disease dynamics. Greek letters represent rates of flow between states over time.}\label{SEIR}
\end{figure}

We have incorporated public health outcomes in the optimization problem objective using a simplified model of disease dynamics, so we should evaluate the efficiency of our vaccine allocation strategy in a more realistic way. We therefore use a more detailed compartmental model that more closely captures the empirically observed disease trends in LAC as a simpler model could lead to inaccuracies in capturing the dynamics \cite{suen2017risk}. Although we employed a simplified model to calculate the herd immunity threshold, we verified that this threshold still holds when applied to our more complex model.. We formulate a compartmental model that includes the following health and treatment states as shown in Figure~\ref{SEIR}: susceptible (S), vaccinated (V), exposed (E), identified infectious (I), unidentified infectious (U), hospitalized (H), recovered (R), and dead (D). Each of these states is stratified by HDs to track variation in disease outcomes by region.

We assume full immunity will be achieved after the first dose of vaccination. We assume that only transmission rates ($\beta$) and vaccination rates ($\xi$) vary geographically by HD. This compartmental model can also be represented as a system of ordinary differential equations, which is shown in Appendix~\ref{compartmental_formulation}.

We use publicly available COVID-19 cumulative cases ($C$), cases ($I$), cumulative deaths ($D$), and hospitalizations ($H$)\cite{covidla} to calibrate the compartmental model. We hereafter use ` $\hat{}$ ' to represent simulated outcomes from the compartmental model to clearly differentiate them from empirical surveillance report values. The calibration objective is to minimize the error in the simulated cumulative cases ($\hat{D}+\hat{H}+\hat{R}+\hat{I}$), cases ($\hat{I}$), cumulative deaths ($\hat{D}$), and hospitalizations ($\hat{H}$). We consider the simulation to start on March 1st, 2020 with a 300-day calibration period. Within the calibration interval, we define five calibration intervals in which calibration parameters are allowed to change between intervals. These five intervals are Lock-down (45 days), social events (45 days), re-opening (45 days), second wave (75 days), and holiday season (90 days). 

We first calibrated the model without vaccinations. We choose the calibrated parameters: $\mathbf{p}=\{\beta_I,\beta_U,\delta_I,\delta_U,\gamma_I,\gamma_U,\gamma_H,\mu_I,\mu_U,\mu_H\}$. These parameters are calibrated using ordinary least squares (OLS) on the following objective:

\begin{align*}
    \min_{\mathbf{p}} & \ \  w_D\sum_{u} \sum_{t} \bar{w}_u(|\hat{D}_{u,t} - D_{u,t}|)+w_C\sum_{u} \sum_{t} \bar{w}_u(|\hat{D}_{u,t}+\hat{H}_{u,t}+\hat{R}_{u,t}+\hat{I}_{u,t} - C_{u,t}|)\\& \ \ + w_I\sum_{u}\sum_{t} \bar{w}_u(|\hat{I}_{u,t}- I_{u,t}|)+w_H\sum_{u} \sum_{t} \bar{w}_u(|\hat{H}_{u,t}- H_{u,t}|)\\
    s.t. & \ \  \text{Compartmental Model System Dynamics in Appendix Section~\ref{compartmental_formulation}}
\end{align*}

In the OLS, four calibration targets are minimized: cumulative death, cumulative identified infections, identified infections over time, and hospitalizations, where all calibration targets are at HD level. $w_d,w_C,w_I,w_H$ are weights for different calibration objectives. The weight $\bar{w}_u$ for each HD $u$ is determined by comparing its total population to the total population in LAC. $\hat{H}_{u,t}, \hat{D}_{u,t}, \hat{R}_{u,t}, \hat{I}_{u,t}$ are the number of hospitalizations, deaths, recovered, and identified infections from the compartmental model for HD $u\in\mathcal{H}$ at time $t$. and $D_{u,t}, C_{u,t}$ (cumulative identified cases for HD $u$ at time $t$), $I_{u,t}, H_{u,t}$ are the actual number of deaths, cumulative identified infections, identified infections, and hospitalizations for HD $u\in\mathcal{H}$ at time $t$, which can be computed from the publicly available data \cite{covidla}.

We calibrate the model sequentially on each of the calibration intervals listed above. The calibrated parameter set from the previous interval is used as the initialization parameter set for the next calibration interval. The calibration results are shown in Appendix~\ref{calibration_results}.

We use this calibrated compartmental model to compute the number of infections averted under different \textit{mega-site} placement and vaccine allocation scenarios (shown in Algorithm~\ref{evaluation}). To do this, we first compute the number of cumulative infections between March and December 2020 with no vaccinations ($N_{CI}$). Then we compute the number of cumulative infections by the end of the year 2020 from the compartmental model with a given vaccination strategy $\pi$ ($\bar{N}_{CI}^{\pi}$), assuming that the vaccination starts 6 months prior to the end of the year 2020 --- after the lockdown has ended and businesses have started to reopen. To evaluate the vaccination program $\pi$'s efficiency in reducing the disease burden, we then calculate the number of infections averted by $N_{CI}-\bar{N}_{CI}^{\pi}$. By running the simulation, the number of cumulative identified infections ($N_{CI}$) is 625,737. We then use this number to compute the number of infections averted.

\begin{algorithm}
\caption{Computing Infections Averted}\label{evaluation}
\begin{algorithmic}[1]
\State Get the optimal solution of $y_{uv}^t,z_{uvw}^t$ from model P\ref{P2_formulation}
\State Compute the vaccination rate $\xi_u^t$ for month $t$ for each HD $u$: $\xi_u = \sum_{v\in \mathcal{H}}y_{uv}^t+\sum_{w\in\mathcal{H}}z_{uvw}^t$ 
\State In the calibrated compartmental model, set $\xi_u(t)$ for each HD $u$ at time $t$, rerun the model, and obtain the number of cumulative infections at the end of year 2020, denoted as $\bar{N}_{CI}$
\State \Return{$N_{CI}-\bar{N}_{CI}$}
\end{algorithmic}
\end{algorithm}

\subsection{Model Results}

Traditional facility location models \cite{marianov2001optimal,silva2008locating,chan2016optimizing,jia2007modeling} generally do not capture commuting patterns in \textit{mega-sites} selection. We therefore first evaluate outcomes using formulation P0, which uses the traditional objective of only considering travel inconvenience assuming all individuals are non-commuters. Next, we solve our model integrating commuting patterns in estimating travel inconvenience (P\ref{P1_formulation}). We then solve a model with both travel inconvenience and public health outcome objectives (P\ref{P2_formulation}) to understand how health and equity objectives change our solution. We assume the importance weight on the public health objective ($\lambda$) is 9, the weight on the smoothness objective ($\lambda'$) is 150 multiplied by the average population across HDs, and the prioritization score ($\delta_u^t$) is $\delta_u^t=0.9\times\delta_u^{t-1}$. These parameters will be varied in sensitivity analyses. We assume the capacity for each \textit{mega-site} is uniformly 400,000 vaccinations per decision period, and each decision period spans a month, with a total of six such periods. 

\subsubsection{Effects of Integrating Commuting Patterns}
In this section, we study how integrating commuting patterns changes the solution of \textit{mega-sites} and vaccine allocation, and how these changes impact travel inconvenience and public health outcomes.

\begin{figure}[h]%
\centering
\includegraphics[width=\textwidth]{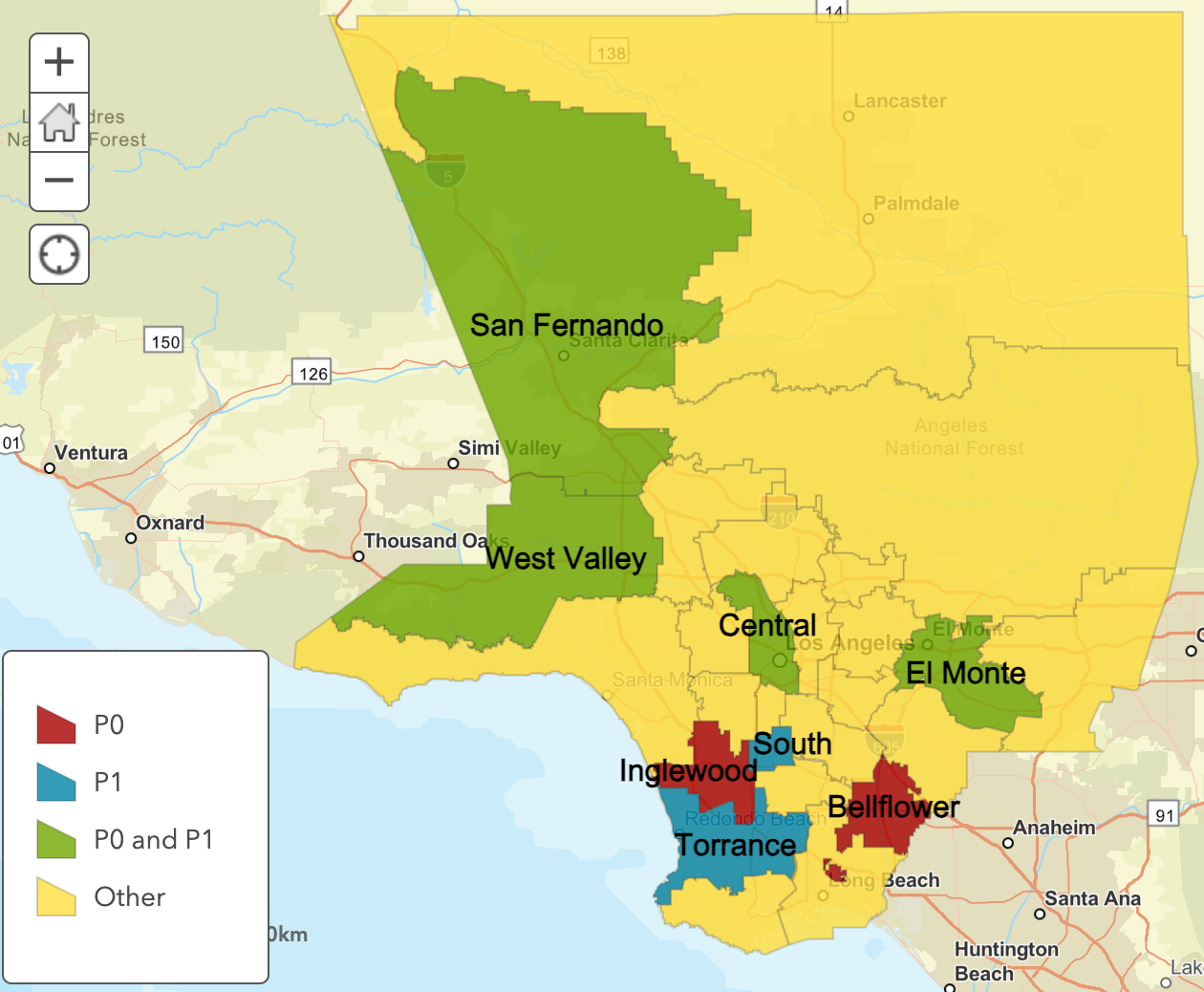}
        \caption{A map of 26 health districts in LAC, with HDs with \textit{mega-sites} selected by P0 marked in red and green, HDs with \textit{mega-sites} selected by P1 marked in blue and green, and HDs without \textit{mega-sites} marked in yellow.}
        \label{fig:locations}
\end{figure}
In the traditional model that minimizes the travel inconvenience only from homes (P0),  San Fernando, West Valley, Central, El Monte, Inglewood, and Bellflower (shown in Figure~\ref{fig:locations}) are selected as locations of \textit{mega-sites}. The model incorporating commuting patterns (P\ref{P1_formulation}) selects four of the same locations as P0 (see Figure~\ref{fig:locations}). The two differences in location (Torrance and South) result in locations that are generally closer to major freeways in LAC (both Torrance and South are closer to I110 and I405 than Bellflower). P\ref{P1_formulation} locations should therefore be more convenient for commuters, evidenced by the substantially reduced total travel inconvenience (nearly 30\% reduction from 276,949,555 minutes to 196,572,213 minutes).

To evaluate and compare public health outcomes between P0 and P\ref{P1_formulation}, we find the number of monthly vaccinations in each HD from the optimization model. We use this value to calculate the HD-specific vaccination rate in the compartmental model, as described in Section~\ref{compart_intro}. Finally, we ran the compartmental to identify the infections averted at the end of 2020 using these vaccination rates (compared to the counterfactual where no vaccination was offered). In the P0 model, 374,538 infections were averted by the end of the year 2020 while the P\ref{P1_formulation} model averted 390,252 infections (4\% more than in the P0 model). This indicates that incorporating commuting patterns in this example also improves disease control even if public health outcomes are not in the objective formulation. For other diseases, particularly those that predominantly affect a minority of the population, further investigation is needed to determine whether incorporating commuting patterns can enhance disease control outcomes.

The inequity objective measures the sum of the maximal difference between HDs with the largest and the smallest vaccinated population percentage during $t\in \mathcal{T}_E$. In this study, we choose $\mathcal{T}_E$ to be the first two months of vaccination due to the high volume of vaccinations observed during the first two months in LAC \cite{covidvaccinedata}. 

Without incorporating inequity into the objective function, we observe large differences in vaccinations across HDs in the first two months. For instance, in the P0 solution, all individuals living in Central are vaccinated in the first month, while none of the individuals living in Pomona are vaccinated. In the second month, there are still zero vaccinations in Pomona, resulting in an inequity score of 2, calculated as $(1-0)+(1-0)$. P\ref{P1_formulation} has a lower inequity score of 1.603 (20\% less than in the P0 model). To further understand how to improve travel inconvenience, public health outcomes, and equity, a model needs to incorporate all of these objectives.

\subsubsection{Multi-Objective Model Results}
In this section, we present results of the optimization formulation (P\ref{P2_formulation}) that integrates objectives in P\ref{P1_formulation} plus public health and equity objectives. We then compare these results with the previously discussed models.

\begin{figure}[h!]
    \centering
    \includegraphics[width=4in]{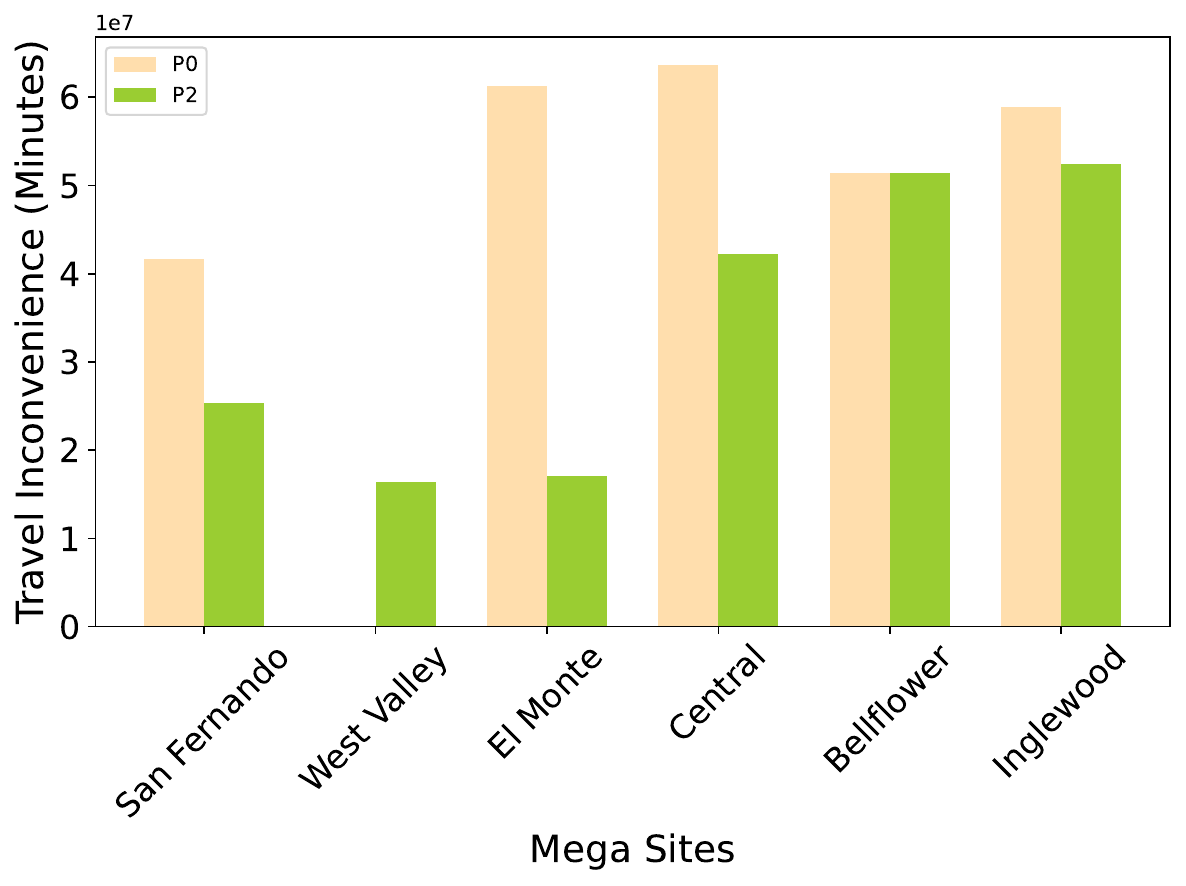}
    \caption{Total travel inconvenience by each \textit{mega-sites} over 6 months. The total travel inconvenience is zero for individuals assigned to West Valley for P0 as it assigns only individuals who live in West Valley to \textit{mega-site} in West Valley.}\label{vaccine_assignment}
\end{figure}

 P\ref{P2_formulation} has the same recommendation as P0 on the locations of \textit{mega-sites}, but the vaccine assignments between the two models are distinct. As shown in Figure~\ref{vaccine_assignment}, the total travel inconvenience in minutes is smaller in P\ref{P2_formulation} for all \textit{mega-sites} except for those receiving vaccination in West Valley. El Monte experiences the most significant decrease, with a reduction of roughly 70 million minutes. In contrast, West Valley faces an increase in travel inconvenience. This is because P0 only assigns individuals who live in West Valley to the \textit{mega-site} in West Valley, which results in zero total travel inconvenience for individuals vaccinated there. On the other hand, P\ref{P1_formulation} assigns individuals who do not live or work near West Valley to West Valley as their more convenient locations are already at capacity, which results in higher total travel time. This increase is offset by decreases in travel inconvenience in other \textit{mega-sites}.

 \begin{table}[h]
\caption{Model results}\label{tab:basecase}%
\begin{tabular}{@{}llll@{}}
\toprule
Model & Travel inconvenience (mins.) & Infections averted &Inequity objective\\
\midrule
P0  & 276,949,555  & 374,538 & 2\\
P\ref{P1_formulation}    & 196,572,213  & 390,252&1.603 \\
P\ref{P2_formulation}  & 204,859,537  & 408,084 &  0.831 \\
\botrule
\end{tabular}
\end{table}

Table~\ref{tab:basecase} shows that formulation P\ref{P2_formulation} is able to achieve improvements in infections averted and lower inequity while still achieving less travel inconvenience than formulation P0 and only a 4\% increase from formulation P\ref{P1_formulation}. The P\ref{P2_formulation} locations are able to prevent approximately 18,000 more infections than P\ref{P1_formulation} by the end of year 2020 and achieve a 48\% reduction in inequity compared to the P\ref{P1_formulation} solution.

Additionally, we tested the model using weekly decision intervals and found that the selected \textit{mega-sites} remained unchanged from the solution of P\ref{P2_formulation}, suggesting that the solution of P\ref{P2_formulation} is not sensitive to a more frequent decision interval within a given range.

Moreover, we conducted a sensitivity analysis assuming immunity is achieved one month after vaccination, considering both one- and two-dose vaccines. We found that the selection of \textit{mega-sites} remains unchanged. Additionally, we considered a scenario where vaccinated individuals have only an average protection of 80\%, reflecting the reality that not all vaccines are 100\% effective. We found that the selection of \textit{mega-sites} remained consistent. These sensitivity analyses suggest that the locations the model recommended are not sensitive to the vaccine immunity threshold within a given range.

\subsubsection{Benchmark Solution: Empirical Locations in LAC}
To benchmark our model against a more realistic setup, we compare model outcomes to those if we fixed the \textit{mega-site} locations to those actually used in LAC from January to May of 2021. More specifically, we solve the optimization problem P\ref{P2_formulation} while fixing the $x_u$ variables to be these empirical \textit{mega-sites} locations. 

In the empirical location model, the \textit{mega-sites} locations are San Fernando, Central, Northeast, Southeast, Inglewood, and San Antonio. Except for one \textit{mega-site} in San Fernando, all other \textit{mega-sites} are located around the Central district. Compared with the empirical \textit{mega-sites} locations, P\ref{P2_formulation} recommends a more dispersed distribution of \textit{mega-sites}, with two \textit{mega-sites} in the northwest, two around the Central district, and two in the east.

The empirical model incurs a total travel inconvenience of 288,467,064 minutes. Compared with the travel inconvenience from the P\ref{P2_formulation} model, the empirical model costs 83,607,527 more minutes in total, or a 40\% increase in total travel time. 87.82\% of individuals are able to be vaccinated at their most convenient \textit{mega-sites}. This represents a 3.1\% decrease from the P\ref{P2_formulation} solution-- approximately 310,000 fewer individuals can now go to the most convenient \textit{mega-sites}.

In the empirical model, 407,855 infections were averted by the end of the year 2020, approximately 200 less than the P\ref{P2_formulation} solution. However, it has an inequity objective value of 0.721, which is less than that in P\ref{P2_formulation}. This indicates there is a trade-off between infections averted and equitable distribution.

Comparing both the results from the public health perspective and transportation perspective, our P\ref{P2_formulation} formulation improves upon the single objective model P\ref{P1_formulation} as well as the model with empirical locations, indicating that both accessibility and public health outcomes can be improved from the empirical setup.

\subsection{Sensitivity Analyses}
Recall the public health objective in our optimization model is minimizing $\lambda \sum_{t\in\mathcal{T}} \sum_{u\in\mathcal{H}} \delta^t_u \zeta_u^t+\lambda'\sum_{t\in\mathcal{T}_E} \tau_t$, which considers prioritization in vaccination allocation ($\delta^t_u$) and equity in distribution for the early phase $\mathcal{T}_E$. In this section, we vary the importance weight on equitable distribution and importance weight on public health objectives to understand how the allocation changes with respect to changes in objectives. We made several assumptions when designing our health objective term, so we tested whether our health objective was reasonable by assessing its correlation with with the number of infections averted. 

\subsubsection{Testing Model Assumptions: Evaluating the Health Objective Term for Concordance with Cases Averted}

Our public health objective consists of vaccination targets (from a simplified SIR model) and the prioritization score (from commuting patterns in LAC). This objective serves as a proxy for true health goals such as infections averted. We evaluate the validity of this approximation in the LAC COVID-19 context in this section. We do so by varying the importance of the health objective term ($\lambda$) and testing if the number of infections reduced is positively correlated to $\lambda$. A good health objective would avert more infections with increasing emphasis on the health objective. This can also serve as a sensitivity analysis on $\lambda$ (relative importance of health objective versus travel inconvenience and inequity score).

\begin{figure}[h!]
     \centering
     \includegraphics[width=0.5\textwidth]{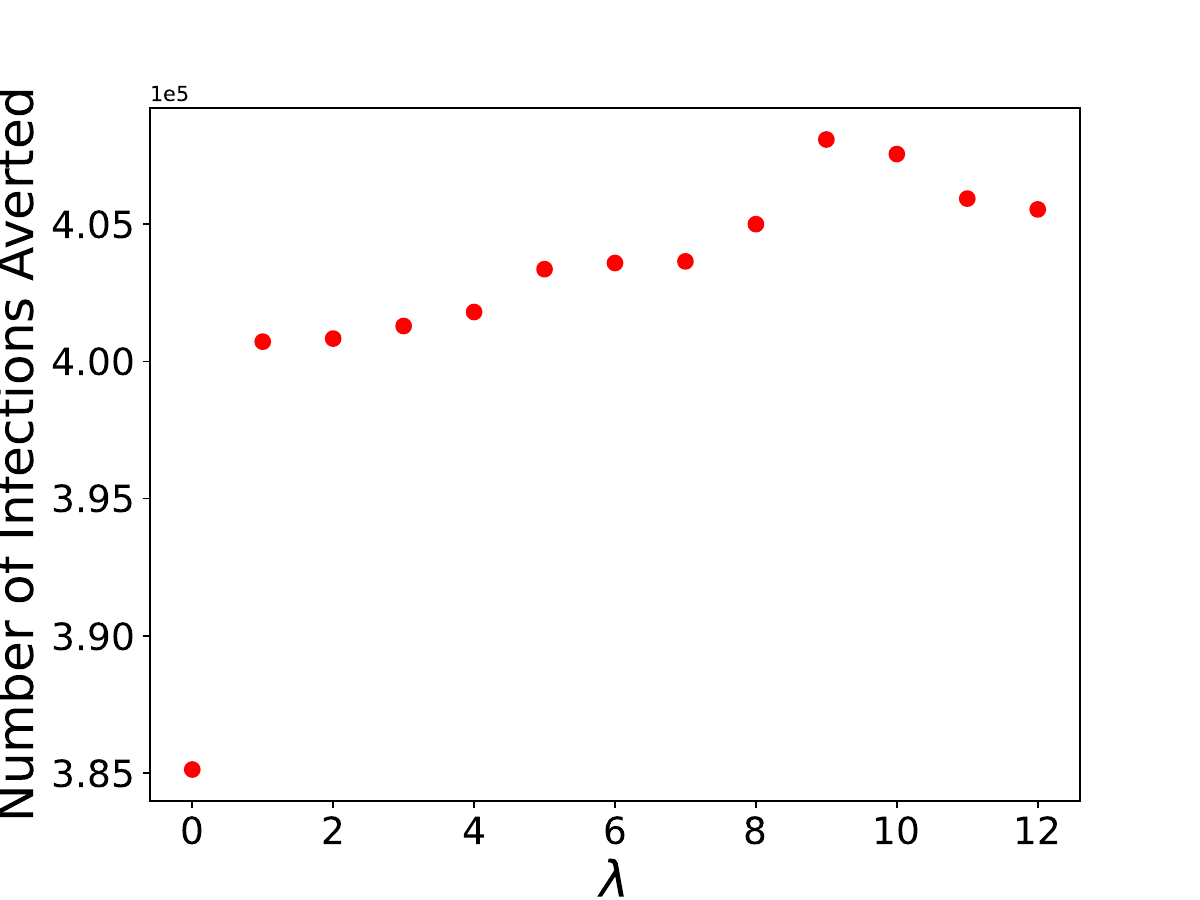}
         \caption{Changes in the number of infections averted with varying $\lambda$.}
         \label{fig:vary_lambda_infections}

\end{figure}

Our health objective term is a reasonable proxy for cases averted when $\lambda\leq 9$. Figure~\ref{fig:vary_lambda_infections} shows infections averted when varying $\lambda$. When $\lambda\leq 9$, the number of infections averted increases with $\lambda$. This suggests that the public health objective we have adopted serves as an effective proxy for the real-world goal of reducing infections. However, when $\lambda$ surpasses 9, the number of infections averted plateaus and even slightly drops, implying that beyond a certain point of weight on the public health objective, the public health objective does not approximate the reduction in infections well. In those cases, when $\lambda$ is exceedingly high, the optimization model prioritized the vaccines only for areas with high risks of becoming infectious while other areas with no vaccine access remain vulnerable, thus not reducing infections effectively.

As expected, when $\lambda$ increases, the travel inconvenience contributes less to the overall objective and the optimal solutions have a larger travel inconvenience. The total travel inconvenience increases by about 4\% when $\lambda$ increases from 0 to 9. Compared with the case when only travel inconvenience is considered in the objective (i.e., $\lambda=0$), the locations of \textit{mega-sites} in southern LAC are more dispersed along the east-west axis. In general, however, the selection of \textit{mega-sites} locations shows limited variation with $\lambda$. 

When $\lambda$ is small, a low inequity score can be achieved. However, the inequity score grows dramatically when $\lambda$ surpasses five. This shows that there exist trade-offs among transportation inconvenience, averted infections, and equity in distribution. We discovered that setting $\lambda$ to 9 yields satisfactory outcomes for all objectives, which validates our choice of $\lambda$ for the P\ref{P2_formulation} model. For future problems, $\lambda$ should be selected within a range from zero to the maximum commuting time per person.


\subsubsection{Changes in Importance Weight on Equitable Distribution}

We use the hyperparameter $\lambda'$ (importance on equitable distribution) to balance between prioritization and equitable distribution. To understand how the importance of equitable distribution influences the solution, we varied $\lambda'$ from 0 to 1000 multiplied by the average population across HDs while other hyperparameters remained the same.

Locations of \textit{mega-sites} do not change when varying $\lambda'$ in the P\ref{P2_formulation} model. However, the inequity score drops from 1.96 to 0 as $\lambda'$ increases from 0 to 250 multiplied by the average population across HDs. By more strongly weighting equitable distribution, the optimization model is configured to assign individuals for vaccination more uniformly across regions within each period.  

\begin{figure}[h!]
     \centering
     \includegraphics[width=0.5\textwidth]{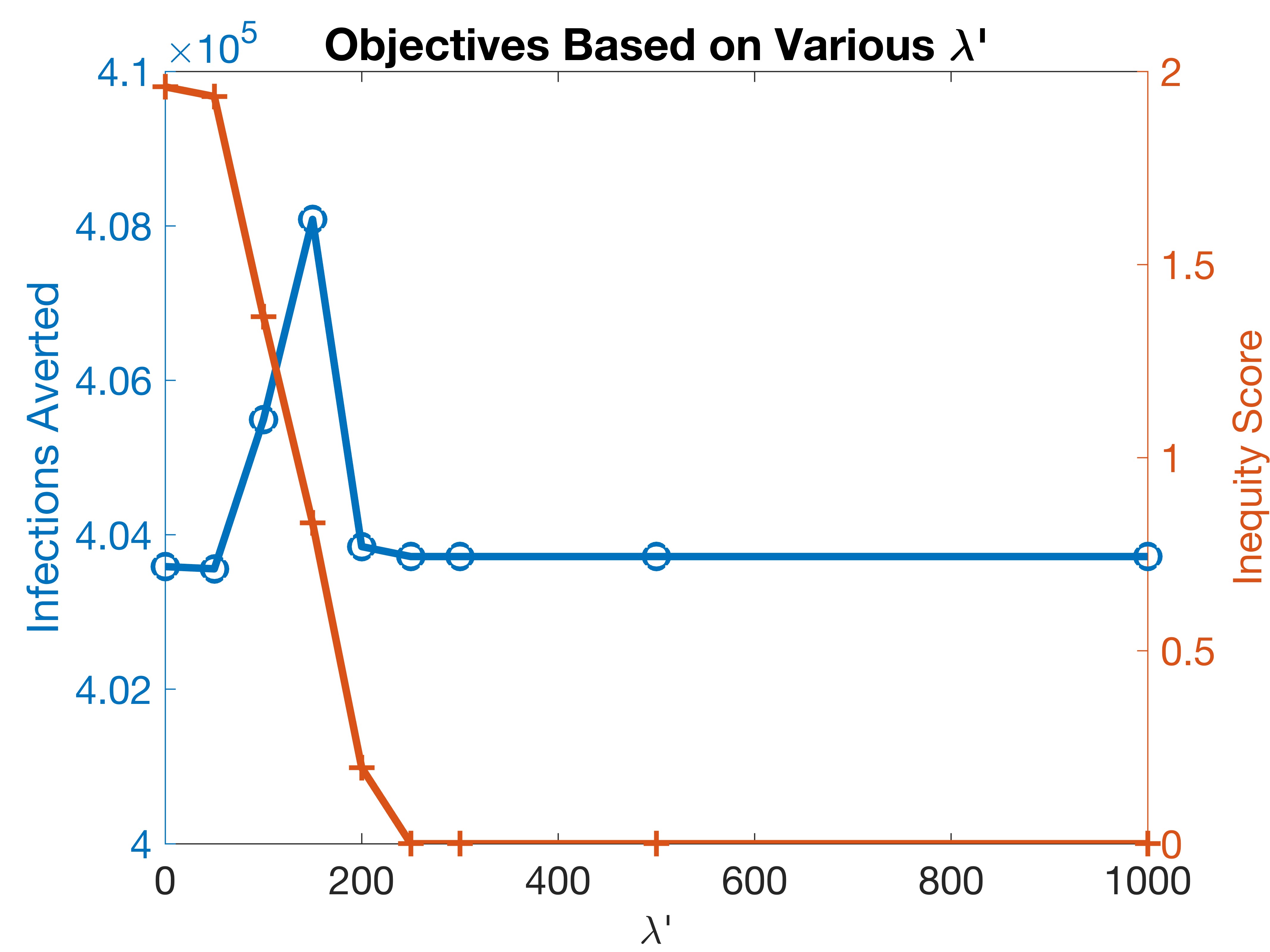}
         \caption{Changes in the number of infections averted and inequity score with varying $\lambda'$.}
         \label{fig:vary_lambda1}

\end{figure}

The number of infections averted is 403,588 when $\lambda'$ is 0, and 403,718 when $\lambda'$ is 1000 multiplied by the average population across HDs, but is higher when $\lambda'$ is 150 multiplied by the average population across HDs (408,084 infections averted). When not considering equitable distribution at all (i.e., $\lambda'=0$), the optimization model prioritized vaccines only to areas at high risk of infection (higher transportation activities) while low-risk areas became more vulnerable due to limited access to the vaccination. The inequity objective value is 1.96 in this case, meaning that the sum of the maximal difference between HDs with the largest and the smallest vaccinated population percentage during the first two months is over 196\%. This indicates a significant inequity in vaccine allocation across HDs. Therefore ignoring equity in the objective not only results in unequal distribution but also reduces the number of averted infections. When the importance of equitable distribution overrides the prioritization for areas at a higher risk of becoming infected (i.e., $\lambda'>150$ multiplied by the average population), the optimization assigns individuals more uniformly for vaccination (thus leads to less than 0.2 inequity objective value), which leads to an ineffective way of distributing limited resources. Combining the results from transportation, public health, and equity perspectives, we choose $\lambda'=150$ multiplied by the average population to balance transportation, prioritization, and equitable distribution across HDs. For future problems, we recommend performing a grid search on $\lambda'$ over the range of 0 to 1000, multiplied by the average population across all regions.

\subsection{Varying the Number of \textit{Mega-Sites}}
We choose the number of \textit{mega-sites} to six to benchmark with empirical solutions. In this section, we vary this number by considering five and seven \textit{mega-sites}. 

When there can be seven \textit{mega-sites}, the selected locations are San Fernando, West Valley, El Monte, Central, West, South, and Bellflower --- only West and South are new and replaced the $mega-site$ in Inglewood. This results in a travel inconvenience of 181,432,015 (11\% less compared to six \textit{mega-sites}), infections averted at 407,505 (0.1\% less compared to six \textit{mega-sites}), and an inequity score of 0.556 (33\% less compared to six \textit{mega-sites}).

When there can be five \textit{mega-sites}, the selected locations are West Valley, El Monte, Central, Bellflower, and Inglewood --- San Fernando is removed, while the other locations remain the same. This leads to a travel inconvenience of 244,839,732 (20\% more compared to six \textit{mega-sites}), infections averted at 387,521 (5\% less compared to six \textit{mega-sites}), and an inequity score of 1.023 (23\% more compared to six \textit{mega-sites}).

Reducing the number of \textit{mega-sites} significantly decreases the infections averted. Increasing the number of \textit{mega-sites} improves travel inconvenience and inequity. However, it does not lead to a notable improvement in the number of infections averted.

\subsection{The Influence of Travel Inconvenience on Vaccine Acceptance}
Prior work indicates that barriers to access vaccination sites, reflected in higher travel inconvenience, may deter individuals from becoming vaccinated \cite{lacombe2018human,schmitzberger2022identifying}. We may therefore be undervaluing the true benefit of reducing travel inconvenience in our models. To remedy this, in this section, we incorporated a likelihood function that computes an individual's probability of choosing to get vaccinated based on their travel inconvenience. We re-evaluate the number of infections averted using the compartmental model after accounting for this additional increase in seeking vaccination. Specifically, in the compartmental model, we compute the vaccination rate ($\xi$) for each month based on the results from the optimization model. Now, if we take into account the willingness to vaccinate as a function of distance, we first calculate the likelihood ($l_{uvw}$) that an individual will choose to get vaccinated, given that they live in HD $u$, work in HD $v$, and are assigned to get vaccinated in HD $w$. Subsequently, the vaccination rate for individuals who live in $u$, work in $v$, and get vaccinated in $w$ is adjusted to $\xi_{uvw} l_{uvw}$.

We adopt the method from \cite{zhu2020personalized,ghafelebashi2023congestion} for computing the likelihood that an individual chooses to become vaccinated as a function of distance traveled. Let $d_{uvw},u\in\mathcal{H},v\in\mathcal{H},w\in\mathcal{H}$ denote the minimum travel inconvenience for an individual who lives in $u$, works in $v$, and assigned to get vaccinated at $w$. We first compute the total difference in the travel inconvenience between becoming vaccinated at $w$ and at another alternative \textit{mega-site} $w'\in\mathcal{H}$:
\begin{equation*}
    \nu_{uvw} = \frac{\sum_{w'\in\mathcal{H}}(d_{uvw'}-d_{uvw})}{\max_{w'}{d_{uvw'}-\min_{w'}{d_{uvw'}}}}
\end{equation*}
where the score is normalized by the difference between the largest and smallest travel inconvenience for this individual to get vaccinated across available \textit{mega-site} locations. Then, we compute the likelihood of an individual choosing to become vaccinated using a logit function:
\begin{equation*}
    l_{uvw}=\frac{exp(v_{uvw})}{1+exp(v_{uvw})}
\end{equation*}
where $l_{uvw}$ denotes the probability of an individual who lives in $u\in\mathcal{H}$, works in $v\in\mathcal{H}$, and chooses to become vaccinated in $w\in\mathcal{H}$. We scale the HD-specific vaccination rates in the compartmental model accordingly to evaluate the number of cases averted.

Given that the decrease in travel inconvenience leads to more people choosing to get vaccinated, 3,813 cases would be averted due to greater access to vaccination using model P\ref{P2_formulation} instead of the empirical setup, which now averts 399,016 infections. This is greater than the unadjusted benefit (an additional 229 cases averted when comparing the unadjusted P\ref{P2_formulation} with the unadjusted empirical solution). Because the empirical solution substantially increases travel inconvenience compared to the P\ref{P2_formulation} solution, fewer individuals tend to come to \textit{mega-sites} for vaccination and thus leads to a worse public health outcome.



In adjusted P\ref{P2_formulation}, the public health outcome in our models outperforms the empirical model. This solidifies the decision to employ a multi-objective formulation that accounts for travel inconvenience, prioritization, and equitable distribution when optimizing \textit{mega-sites} and allocations. 

\section{Conclusion}\label{conclusion}

We present an optimization formulation for selecting locations of vaccine sites during the early stage of the vaccination. In the optimization model, we integrate multiple objectives including travel inconvenience, prioritization, and equitable distribution. 

In formulating the transportation objective, we consider commuter patterns within LAC. We use the estimated origin-destination pair to infer the number of commuters and non-commuters in each HD. We compute the travel inconvenience for the vaccination for each origin-destination pair, considering the accessibility of vaccination sites from both residential and workplace perspectives.

When modeling the public health objective, we establish a straightforward linear goal to pinpoint and prioritize high infection risk zones. This approach guarantees that the optimization model remains solvable and manageable. We also factor in the fair allocation of vaccines across HDs during the initial phases, resulting in a vaccination strategy that upholds equity.

Our numeric case study on COVID-19 vaccination in LAC demonstrates that our model significantly eases travel burden and enhances public health outcomes compared to the empirical approach. It cuts down the total transportation time for vaccinations by roughly 26\% for LAC residents. Furthermore, we contrast the empirical solution's simulation results on public health outcomes with those from our optimized strategy; our model shows a reduction in infections by an additional 229 cases. When accounting for the influence of travel inconvenience on vaccine acceptance, our model demonstrates a reduction of an additional 3,813 infection cases.

Furthermore, we conduct sensitivity analyses on the significance of equitable distribution weight, and the impact of public health outcome weight. We demonstrate how travel inconvenience and the number of infections averted changes with respect to changes in these parameters. Additionally, our results reveal that the public health objective value in our model mirrors the variations in the number of infections prevented. We also factor in each individual's propensity to come to the vaccine site for vaccination, based on their designated location and commuting behaviors. The outcomes continue to indicate that our model surpasses the empirical solution in the effectiveness of reducing infections.

The policy implications for our work are (1) locations of \textit{mega-sites} for dispensing a large number of vaccinations are recommended to be placed widely across the region for reducing the overall travel inconvenience, and (2) merely minimizing travel inconvenience in determining vaccine site locations and vaccine distribution falls short -- it's also essential to identify and prioritize high infection risk areas to more effectively control disease transmission (3) balancing travel inconvenience, the number of infections prevented, and equity is challenging. However, with an appropriate setup in the optimization model, it's possible to achieve satisfactory outcomes in all these aspects.

We must acknowledge several limitations of this work. While in reality, vaccines were not available until later, using July 2020 in our example allows us to study what would occur without the changing dynamics due to the real-world implementation of vaccine uptake (we can calibrate the model to non-vaccine data) but is close enough to the actual vaccine available date that we can compare between our optimal policy and real-world policies. We assume the capacities of different vaccination sites are uniform. In reality, different places might have different capacities based on their size and number of nurses. Future research could focus on addressing the problem of optimizing capacity, which may vary across distribution sites. In practice, capacity constraints were most significant at the very beginning of when the vaccine became available, while there will be more supply over time as local clinics begin to offer the vaccine. We only consider individuals' residence/work HD in vaccine assignment. However, additional factors like age and work industry could also be taken into consideration. It is challenging to accommodate all commuters during peak morning and evening hours at commuter-dominated \textit{mega-sites}, requiring some to schedule vaccinations outside these preferred times. Future research can incorporate daily scheduling into the model.

Our model assumes that all individuals receive their vaccinations at \textit{mega-sites}. In practical terms, these \textit{mega-sites} are operational primarily in the initial stages and are eventually phased out as vaccinations become available at local clinics. Nonetheless, this is unlikely to significantly affect our findings since historically, the majority of vaccinations occur at \textit{mega-sites}, and minor variations in the number of people getting vaccinated through the optimization model do not substantially alter the solution. Our model incorporates varying levels of detail to reflect differences in the availability of information, data, or preferences across its components. The vaccination rate in our model is higher but still comparable to empirical data. We assumed full vaccination in the model, which is not realistically achievable. While travel time within larger HDs like San Fernando was not initially considered, a sensitivity analysis incorporating these factors demonstrated that the solution held consistent. We acknowledge that the herd immunity level is subject to uncertainty due to assumptions made and uncertainty in disease dynamics. However, we performed additional analyses using 80\% and 100\% herd immunity thresholds and found that the selection of \textit{mega-site} locations remained stable. Additionally, herd immunity thresholds tend to be similar across various diseases; for example, in our COVID-19 case study, we calculated an average herd immunity threshold close to 97\%, similar with that of measles \citep{ashby2021herd}. 

Despite these limitations, we believe that this work provides interesting insights into not only the locations of vaccine sites but also the allocation of vaccines. Our paper provides a tractably solvable optimization formulation and a numerical example to show how to reduce the travel burden and prevent more infections from the empirical approach. These results provide insight into future work on improving vaccination strategies to manage the disease burden in scenarios where an outbreak is ongoing and a new vaccine has just been introduced. While prior work on optimal location focused primarily on inconvenience, our model incorporates health considerations (such as prioritization) and commuting patterns, adding realism to the traditional transportation problem.

\backmatter





\bmhead{Acknowledgements}

This material is based upon work supported by the National Science Foundation under Grant No. 2237959. Research reported in this publication was also supported by the National Library of Medicine of the National Institutes of Health under award number 5R21LM013697. The content is solely the responsibility of the authors and does not necessarily represent the official views of the National Institutes of Health. Research reported in this publication was also supported by the National Research and Development Agency of Chile under award number AFB230002.





\newpage
\begin{appendices}

\section{Compartmental model dynamics}\label{compartmental_formulation}

\begin{align*}
    &\frac{dS_{v}}{dt}=- \sum_u {\beta_I}_{uv} (t)\times S_v(t) \times \frac{{I}_u(t)}{{N}_u(t)}- \sum_u {\beta_U}_{uv}(t) \times S_v(t) \times \frac{U_u(t)}{{N}_u(t)} - \xi_v (t) \times S_v(t)\\
    &\frac{dE_v}{dt}= \sum_u {\beta_I}_{uv}(t) \times S_v(t) \times \frac{{I}_u(t)}{{N}_u(t)}+ \sum_u {\beta_U}_{uv}(t) \times S_v(t) \times \frac{U_u(t)}{{N}_u(t)}-({\delta_I}(t)+{\delta_U}(t)) \times E_v(t)\\
    &\frac{dI_v}{dt} = {\delta_I(t)} \times E_v(t)-({\gamma_I(t)}+{\mu_I(t)}+\eta(t)) \times I_v(t)\\
    &\frac{dU_v}{dt} = {\delta_U}(t) \times E_v(t)-({\gamma_U}(t)+{\mu_U}(t)+\theta(t)) \times U_v(t)\\
    &\frac{dH_v}{dt} = {\eta}(t) \times I_v(t)+\theta(t)\times U_v(t)-({\gamma_H}(t)+{\mu_H}(t)) \times H_v(t)\\   
    &\frac{dR_v}{dt} = {\gamma_I}(t) \times I_v(t)+{\gamma_U}(t) \times U_v(t)+{\gamma_H}(t) \times H_v(t)\\
    &\frac{dD_v}{dt} = {\mu_I}(t) \times I_v(t)+{\mu_U}(t) \times U_v(t)+{\mu_H}(t) \times H_v(t)\\
    &\frac{dV_v}{dt} =\xi_v (t) \times S_v(t)
\end{align*}
The presented compartmental model captures the dynamics of a partially vaccinated population in the context of infectious disease transmission. The population is divided into various compartments: \(S_v(t)\) (susceptible individuals), \(E_v(t)\) (exposed individuals who are infected but not yet infectious), \(I_v(t)\) (identified infected individuals), \(U_v(t)\) (unidentified infected individuals), \(H_v(t)\) (hospitalized individuals), \(R_v(t)\) (recovered individuals), and \(D_v(t)\) (dead individuals). Additionally, \(V_v(t)\) represents the cumulative number of individuals who have been vaccinated. The model tracks the rates of change for each compartment over time (\(t\)) using a set of differential equations. Transmission between individuals is governed by infection rates (\(\beta_{Iuv}(t)\) and \(\beta_{Uuv}(t)\)) based on contact with both identified (\(I_u(t)\)) and unidentified (\(U_u(t)\)) infected individuals in the unvaccinated population. Individuals move between compartments due to various processes such as infection, recovery, hospitalization, or death, characterized by rates such as \(\delta_I(t)\) (rate of exposed individuals becoming identified infected), \(\delta_U(t)\) (rate of exposed individuals becoming unidentified infected), \(\gamma_I(t), \gamma_U(t), \gamma_H(t)\) (recovery rates), and \(\mu_I(t), \mu_U(t), \mu_H(t)\) (mortality rates). Vaccination dynamics are modeled with a term \(\xi_v(t)\), which reduces the number of susceptible individuals by transitioning them into the vaccinated state. 

\section{COVID-19 Model Calibration Results}\label{calibration_results}
\begin{figure}[h!]
\centering
\includegraphics[width=0.4\textwidth]{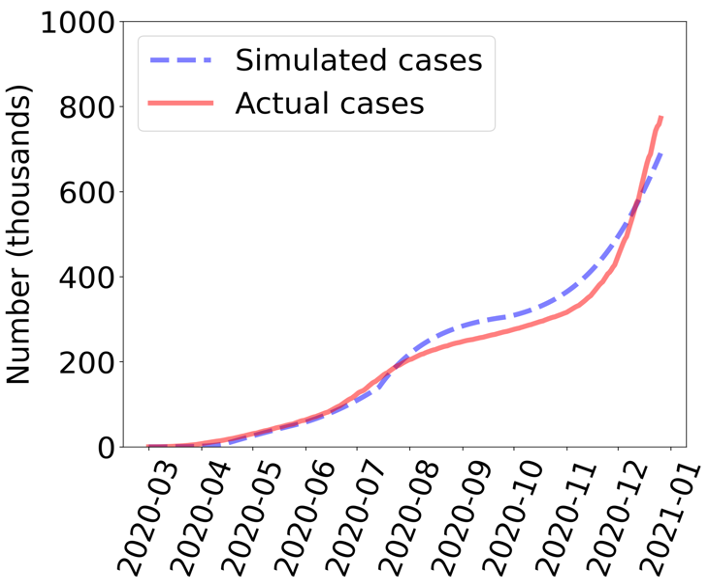}
\caption{Calibration result of cumulative identified infections on the aggregate level}\label{fig:infections_calibration}
\end{figure}
\begin{figure}[h!]
\centering
\includegraphics[width=0.4\textwidth]{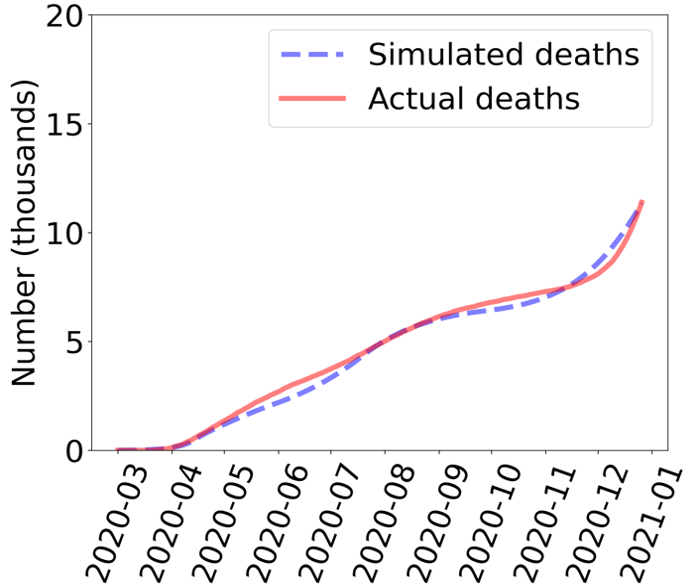}
\caption{Calibration result of cumulative deaths on the aggregate level}\label{fig:deaths_calibration}
\end{figure}
\begin{figure}[h!]
\centering
\includegraphics[width=0.4\textwidth]{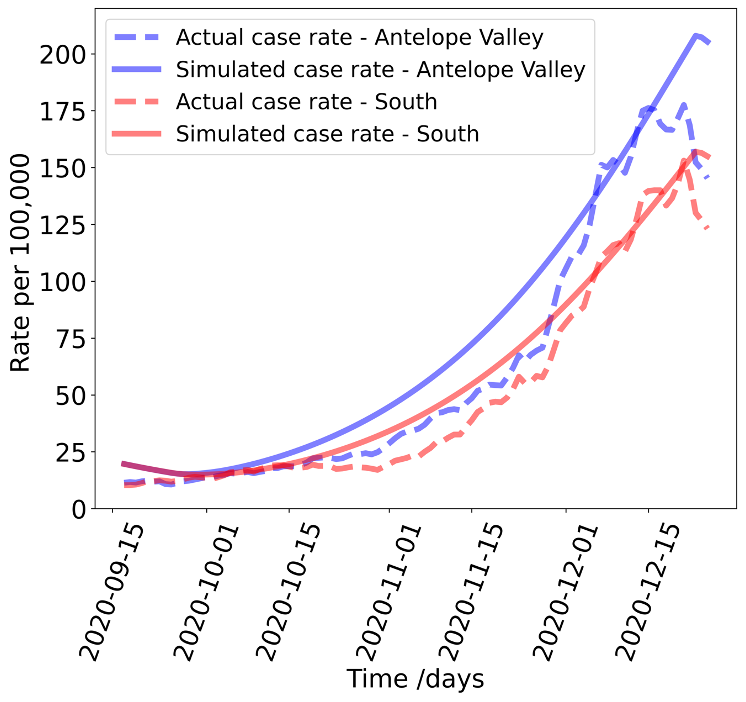}
\caption{Calibration result of case rate on the HD level}\label{fig:caserate_calibration}
\end{figure}
Figure~\ref{fig:infections_calibration} shows the calibration result of cumulative identified infections on the aggregate level. Figure~\ref{fig:deaths_calibration} shows the calibration result of cumulative deaths on the aggregate level. Figure~\ref{fig:caserate_calibration} shows the calibration result of the case rate on the HD level. All calibration results demonstrate that the calibrated compartmental model is able to capture the disease dynamics in both aggregate and HD levels.
\clearpage
\section{Herd Immunity Threshold by HD}\label{herd_results}

\begin{figure}[hbt!]%
\centering
\includegraphics[width=0.6\textwidth]{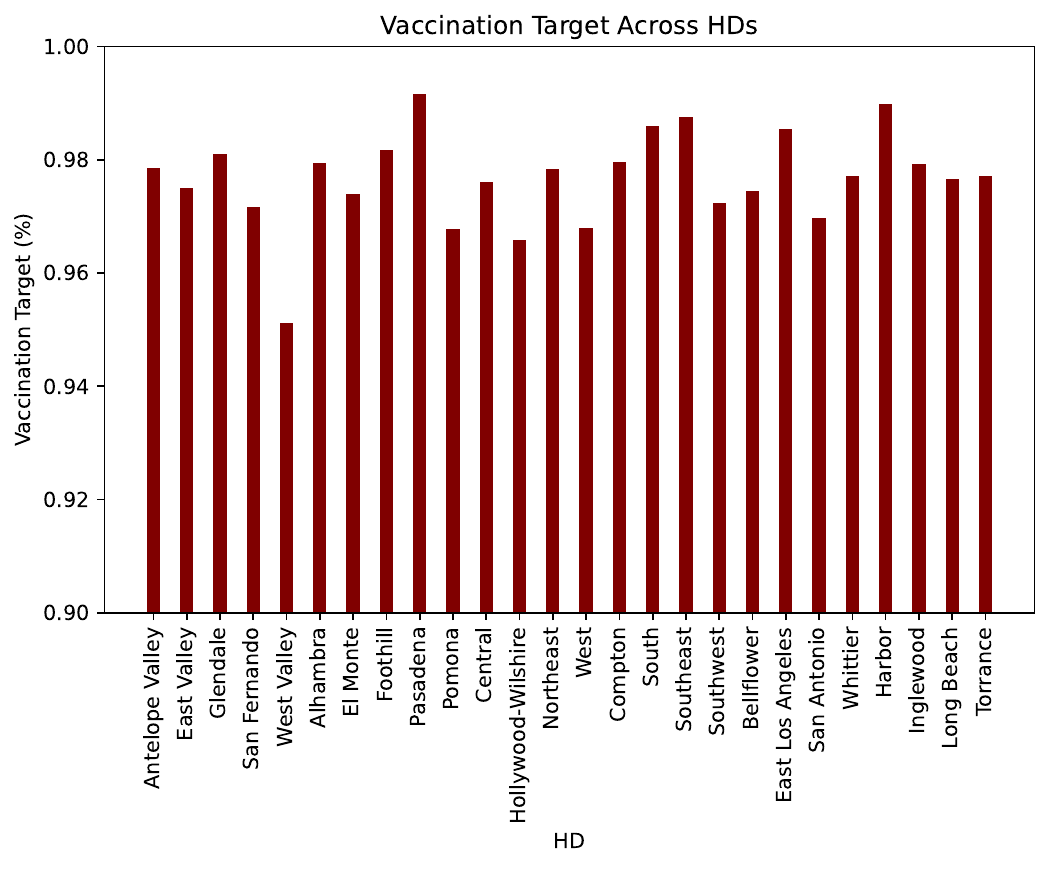}
\caption{Vaccination targets across HDs. }\label{vac_target}
\end{figure}
\section{Choices of Prioritization Scores}

\begin{table}[h]
\caption{Model results}
\begin{tabular}{@{}llll@{}}
\toprule
Model & Travel inconvenience (mins.) & Infections averted &Inequity objective\\
\midrule
Betweenness Centrality  & 210,216,624  & 405,247 & 1.799\\
Population     & 200,576,766  & 404,068&0.590 \\
Case rate  & 201,865,615  & 401,503 &  0.447 \\
\botrule
\end{tabular}
\end{table}



\end{appendices}


\bibliography{sn-bibliography}

\end{document}